\documentclass[11pt, reqno]{amsart}

 \usepackage{amsmath}
 \usepackage{amssymb}

\DeclareMathAccent{\mathring}{\mathalpha}{operators}{"17}

 \usepackage{color}

\newtheorem{theorem}{Theorem}%[section]
\newtheorem{lemma}[theorem]{Lemma}     

\newtheorem{corollary}[theorem]{Corollary}
\theoremstyle{remark}
\newtheorem{remark}[theorem]{Remark}

\theoremstyle{definition}

\newcommand\cbrk{\text{$]$\kern-.15em$]$}}
\newcommand\opar{\text{\raise.2ex\hbox{${\scriptstyle | }$}\kern-.34em$($} }

\newcommand{\Vol}{\text{\rm Vol}\,}

 \makeatletter
 \def\dashint{%  
 \operatorname%
 {\,\,\text{\bf--}\kern-.98em\DOTSI\intop\ilimits@\!\!}}
\def\ninf{\qopname\relax\@empty{\!\!\phantom{p}\!\!\!inf}}
 \makeatother

\newcommand\bbeta{\text{\raise-.2ex\hbox{$\bm{\beta}$}}}

 \makeatletter
 \def\dashint{%  
 \operatorname%
 {\,\,\text{\bf--}\kern-.98em\DOTSI\intop\ilimits@\!\!}}
\def\ninf{\qopname\relax\@empty{inf\phantom{p}\!\!\!}}
 \makeatother

\newcommand\bR{\mathbb{R}}

\begin{document}

\title[Brouwer's fixed point theorem]
{On a new proof of the key step in the proof of Brouwer's fixed point theorem}

\author{N.V. Krylov}
 
\email{nkrylov@umn.edu}
\address{127 Vincent Hall, University of Minnesota,
 Minneapolis, MN, 55455}



\begin{abstract}
We present a solution of   Exercise 1.2.1 of \cite{Kr_08}
which yields a short new proof of a key step
in one of proofs of
Brouwer's fixed point theorem, 1910.
A few people asked the author about the details
of the solution and they might be interesting to a broader audience.
Our approach is absolutely different from
the ones using algebraic or differential topology or differential calculus and is
based on a simple observation which somehow
escaped many authors treating this theorem
in the past.

\end{abstract}

\maketitle

By $\bR^{d}$ we denote the Euclidean
space of points $x=(x^{1},...,x^{d})$. When it makes sense,
for real-valued $u(x)$ on $\bR^{d}$ we denote 
$$
D_{i}u =\frac{\partial u}{\partial x^{i}}.
$$

If $F=(F^{i})$ is a smooth mapping of $\bR^{d}$
to $\bR^{d}$, we set
$$
DF=(a^{ij} )_{i,j=1 }^{ d},\quad a^{ij}=D_{j}F^{i}.
$$

Here is a key   lemma in one of the the proofs 
of
Brouwer's fixed point theorem, see, for instance,
pages 467--470 in
N. Dunford and J.T. Schwartz \cite{DS_58} or  W. Kulpa \cite{Ku_89},
in which \cite{DS_58} is not referred to. In both references $F$ is assumed to be in $C^{2}$
since they use integrating by parts unlike us. We present
the proof of the lemma with  some intentional gaps, closed
later in   Remark \ref{remark 3.21.2}.
The easy way, Brouwer's fixed point theorem is derived from
the lemma, is well known and can be found in
\cite{DS_58}, \cite{Ku_89}, \cite{Kr_23} and some other sources.
\begin{lemma}
                           \label{lemma 3.21.1}
Let $\Omega$ be a (connected) bounded domain
in $\bR^{d}$
with $C^{1}$ boundary
and let $F, G:\bar\Omega\to\bR^{d}$ be  $C^{1}(\bar{\Omega} )$
mappings such that 
 $$
F=G\quad\text{on}\quad \partial \Omega.
 $$ 
Then
 $$
\int_{\Omega}\det D F\,dx=
\int_{\Omega}\det DG\,dx.
 $$

\end{lemma}

Proof.
Observe that
 for small $t $ the mappings  
$F_{t}=tF(x)+ x $ and $G_{t}=tG(x)+ x$ are one-to-one
on $\bar\Omega$ and,  by the implicit function theorem,
have $C^{1}$-inverse mappings on $\Omega$.
Because of that they
map $\partial \Omega$ onto the boundary of
$F_{t}(\Omega)$ which is $F_{t}(\partial\Omega)$ and is, of course, the same
as the boundary of $G_{t}(\Omega)$. 
Furthermore, for small $t$ the intersection of
$F_{t}(\Omega)$ and $G_{t}(\Omega)$ is obviously nonempty
and since they are connected and have the same boundary,
$F_{t}(\Omega)=G_{t}(\Omega)$ and
 $$
\text{Vol}\,F_{t}(\Omega)=\text{Vol}\,G_{t}(\Omega)
 $$
for small $t$.
We express this equality in terms of $D F_{t}$ and 
$DG_{t}$ as
\begin{equation}
                                      \label{3.21.1}
\int_{\Omega}\det (tD F+I)\,dx=
\int_{\Omega}\det (tDG+I)\,dx,
\end{equation}
 where $I$ is the unit $d\times d$-matrix,
 then use the fact 
both parts of \eqref{3.21.1} 
are polynomials in $t$. Since they coincide
for small $t$, they are identical
 and by comparing the
coefficients of $t^{d}$ we get the result.
The lemma is proved. \qed

Now we   comment on some steps some people
may regard as missing in the proof
of Lemma \ref{lemma 3.21.1}.
\begin{remark}
                            \label{remark 3.21.2}
First question: Why (for small $t$) is $F_{t}$
one-to-one? Assume that there are two
points $x',x''\in\bar\Omega$ such that
$F_{t}(x')=F_{t}(x'')$. Then
by denoting by $N_{F}$ the Lipschitz constant of $F$
and taking $t$ such that $tN_{F}\leq 1/2$, we get
$$
|x'-x''|=t|F_{t}(x')-F_{t}(x'')|\leq N_{F}t|x'-x''|
\leq(1/2)|x'-x''|
$$
and $|x'-x''|=0$.

Second question: Why (for small $t$) does $F_{t}$
map  $\partial\Omega$ onto $\partial F_{t}(\Omega)$?
Assume that there is a point $x_{0}\in\partial\Omega$
such that $F_{t}(x_{0})\not\in\partial F_{t}(\Omega)$.
Then $F_{t}(x_{0})=:y_{0}\in F_{t}(\Omega)$
(which is an open set in $\bR^{d}$
by the implicit function theorem)
and consequently there is $x_{1}\in\Omega$
such that $F_{t}(x_{1})=:y_{0}$. We have 
$x_{1}\in\Omega$, $x_{0}\in\partial\Omega$,
so that $x_{1}\ne x_{0}$, contradicting the one-to-one
property. Thus, $F_{t}(\partial\Omega)\subset
\partial F_{t}(\Omega)$.

That, conversely, any $y_{0}\in \partial F_{t}(\Omega)$ is in 
$F_{t}(\partial\Omega)$ follows from the fact that
there is a sequence $x_{n}\in\Omega$ such that
$y_{n}:=F_{t}(x_{n})\to y_{0}$ as $n\to\infty$
and for any subsequence of $x_{n}$
converging, say to $x_{0}\in\bar\Omega$ we have $F_{t}(x_{0})
=y_{0}$, which leaves only one possibility for $x_{0}$:
$x_{0}\in \partial \Omega$, since $y_{0}\not
\in F_{t}(\Omega)$.

Third question: Why (for small $t$)
does the equality $\partial F_{t}(\Omega)=
\partial G_{t}(\Omega)$ imply that 
$  F_{t}(\Omega)=
 G_{t}(\Omega)$? Here we use that $\Omega$
is connected and first prove that
$$ 
F_{t}(\Omega)\cap
 G_{t}(\Omega)\ne\emptyset.
$$
  
For that we fix any $x_{0}\in\Omega$ and show that
if $t$ is sufficiently small, then 
\begin{equation}
                                        \label{3.21.5}
F_{t}(x_{0})
\in G_{t}(\Omega).
\end{equation}
Let $t$ be so small that, for any $x\in\Omega$,
 the distance
of $x_{0}+tF(x_{0})-tG(x)$ to $\partial\Omega$
is at least half the distance of $x_{0}$ to
$\partial\Omega$. Then define  
$$
x_{k+1}=x_{0}+tF(x_{0})-tG(x_{k}),\quad k=0, 1,...
$$
Decreasing $t$ if necessary we may assume that
$t|G(x')-G(x'')|\leq(1/2)|x'-x''|$ and then it follows that
the sequence $x_{k}$ converges and the limit point,
say $x'$ is in $\Omega$ and satisfies 
$F_{t}(x_{0})=G_{t}(x')$. This proves \eqref{3.21.5}.

Now assume that $F_{t}(\Omega)\not\subset
 G_{t}(\Omega)$. Then there is $y_{1}=F(x_{1})\in 
F_{t}(\Omega)$ such that $y_{1}\not\in 
G_{t}(\Omega)$. Take a broken line $x_{s}$, $s\in[0,1]$,
 inside
$\Omega$ connecting $x_{0}$ and $x_{1}$.
On the one end of this broken line $F(x_{1})
\not\in 
G_{t}(\Omega)$ and on the other $F_{t}(x_{0})
\in G_{t}(\Omega)$. Hence there is $s\in[0,1]$
such that $x_{s}\in \Omega$, $F_{t}(x_{s})
\in F_{t}(\Omega)$ and $F_{t}(x_{s})\in \partial
G_{t}(\Omega)=\partial
F_{t}(\Omega)$, the latter contradicting
$F_{t}(x_{s})
\in F_{t}(\Omega)$. It follows that 
$F_{t}(\Omega)\subset G_{t}(\Omega)$.
By symmetry, $G_{t}(\Omega)\subset F_{t}(\Omega)$
(if $t$ is small enough) and 
$G_{t}(\Omega)= F_{t}(\Omega)$.
\end{remark}
 
 \begin{remark}
An analytic proof of Lemma \ref{lemma 3.21.1} can be obtained
if one proves as in \cite{DS_58} that for smoother $\bR^{d}$-valued $H$
on $\bR^{d}$
$$
\det D H=(1/d){\rm div}\,\hat H,\quad \hat H_{j}=H^{i}A_{ij},
\quad\sum_{j=1}^{d}D_{j}A_{ij}=0,
$$
where $A_{ij}$ are the cofactors of
 $ D_{j} H^{i}  $ in the matrix 
$D H$.

Indeed, in that case
$$
\frac{d}{dt}\det D [tF+(1-t)G]
=\sum_{i,j=1}^{d}A_{ij}D_{j}[F^{i}-G^{i}] 
=\sum_{i,j=1}^{d}D_{j}\big(A_{ij}[F^{i}-G^{i}] \big),
$$
the integral over $\Omega$ of the last divergence
is reduced to the integral over its boundary 
and $A_{ij}[F^{i}-G^{i}]=0$ on $\partial \Omega$.

The above  proof  for smoother
$F,G$ is, actually, based on the Stokes theorem
from the theory of differential forms.

\end{remark}

 For completeness we provide the rest of the
proof of Brouwer's fixed point theorem
following 
pages 467--470 of
N. Dunford and J.T. Schwartz \cite{DS_58}.

\begin{corollary}
                                  \label{corollary 3.21.1}
For the domain $\Omega$ from Lemma \ref{lemma 3.21.1}   
 there is no $C^{1}(\bar{\Omega})$ function
$G:\bar{\Omega} \to\partial \Omega$ such that 
$G(x)=x$ on $\partial \Omega$.

\end{corollary}

Indeed, if we assume the contrary, then for $F(x)=x$
we find
$$
\Vol \Omega=\int_{\Omega}\det DG\,dx.
$$
However,   the condition
$G:\bar{\Omega} \to\partial \Omega$ implies that
all partial derivatives of $G$ are tangent to
$\partial\Omega$, in particular,
all $d$-columns of $DG(x)$ are
tangent to $\partial\Omega$ at the point $G(x)$.
But the tangent plane to $\partial\Omega$
at any point is only $(d-1)$-dimensional,
so that the columns of $DG$ are linearly dependent 
and, hence, $\det DG=0$. This yields
a contradiction, $\Vol \Omega=0$, and proves our claim.\qed

\begin{theorem} 
                                  \label{theorem 3.21.1}
Let $B$ be the closed unit ball
centered at the origin and let $f:\bar{B} \to \bar{B} $ be a 
continuous mapping.
Then
$f$ has fixed points in $\bar{B} $.
\end{theorem}

Proof. First assume that $f$ is smooth.
Suppose that there are no fixed points
 and for each $x\in\bar{B} $
define $G(x)\in \bar{B} $ as 
$$
G(x)=x-t(x)(f(x)-x),
$$
 where
$t(x)\geq0$ is the root of the equation 
$$
|x-t(f(x)-x)|=1.
$$ 
From the geometric picture it is clear
 that this equation has always two distinct roots
($f(x)\ne x$) one is strictly negative and the other
is nonnegative (zero if $x\in\partial B $). This means that 
the discriminant of the quadratic
equation
$$
|x|^{2}-2t(x,f(x)-x)+t^{2}|f(x)-x|^{2}=1
$$
is strictly positive, smooth, and its square root is smooth,
so that
$$
t(x)=\frac{(x,f(x)-x)+\sqrt{(x,f(x)-x)^{2}+(1-|x|^{2})|f(x)-x|^{2}}} 
{|f(x)-x|^{2}}
$$
is a smooth function along with $f(x)$
and $G(x)$.  Then we can use  Corollary
\ref{corollary 3.21.1} and finish the proof in the case
of smooth $f$.

In the general case, let $f_{n}$ be a sequence of polynomials
such that 
$$
|f_{n}- f|\leq 1/n
$$
   on $\bar B$.
By replacing $f_{n}$ with $(n/(n+1))f_{n}$, if necessary,
we may assume that the $f_{n}$'s map $\bar B$ into itself.
Then there exist $x_{n}\in\bar B$ such that $f_{n}(x_{n})
=x_{n}$. Obviously any converging subsequence of $x_{n}$
converges to a fixed point of $f$.
The theorem is proved. \qed

\end{document}